\documentstyle[11pt]{article}

\begin{document}
\renewcommand{\baselinestretch}{2}
\input{psfig}
\large
\centerline {\Large \bf Expanding Direction of the Period Doubling Operator}

\vskip5pt
\centerline{ Yunping Jiang}

\centerline{The Institute for Mathematical Sciences}

\centerline{SUNY at Stony Brook, Stony Brook, 11794, USA}

\vskip5pt
\centerline{ Takehiko Morita}

\centerline{Department of Mathematics, Tokyo Institute of Technology}
\centerline{Ohokayamo, Meguro, Tokyo 152, Japan}

\vskip5pt
\centerline{Dennis Sullivan}

\centerline {The Department of Mathematics, Graduate Center of CUNY}

\centerline {33 West 42nd Street, New York, N.Y. 10036, USA}

\vskip10pt
\centerline{October 15, 1990}

\vskip25pt

\centerline{\bf Abstract}

We prove that the period doubling operator has an expanding direction
at the fixed point. We use the induced operator, a
``Perron-Frobenius type
operator'',
to study the linearization of the period doubling operator at its fixed
point.
We then use a sequence of
linear operators with finite ranks to study this induced operator.
The proof is constructive. One can calculate the expanding direction and the rate
of expansion of the period doubling operator at the fixed point.

\vskip25pt
\centerline{ \bf Contents }

\vskip5pt
\S 1 Introduction.

\vskip5pt
\S 2 The Period Doubling Operator and the Induced Operator.

\vskip3pt

$\hskip10pt$ \S 2.1 From the period doubling operator to the induced
operator.

\vskip3pt

$\hskip10pt$ \S 2.2 The induced operator ${\cal L}_{\varphi}$.

\vskip3pt

$\hskip10pt$ \S 2.3 A general theorem for operators like the
induced operator.

\vskip5pt
\S 3 The Construction of the Expanding Direction.

\vskip3pt

$\hskip10pt$ \S 3.1 An easy observation.

\vskip3pt

$\hskip10pt$ \S 3.2 The construction.

\vskip3pt

$\hskip10pt$ \S 3.3 A program.

\pagebreak

\centerline{\Large \bf \S 1 Introduction.}

\vskip5pt
{\bf Perron-Frobenius operator.}
Suppose $M^{n}$ is a compact, connected,
oriented and smooth Riemannian manifold, $\Omega \subset M^{n}$ is an open
set and  $\sigma :\Omega \rightarrow M^{n} $ is an expanding mapping.
Let ${\cal B}=\{ v|v $ is a
complex Lipschitz vector field on $\Omega \}$ and let $\phi $ be
a real Lipschitz function on $\Omega $.  The Perron-Frobenius operator
${\cal L}_{\phi } : {\cal B} \rightarrow {\cal B}$ is defined by
\[ ({\cal L}_{\phi }v)(x) = \sum_{y\in \sigma^{-1}(x)} (e^{\phi
(y)})v(y)\]
for $v\in {\cal B}$. An eigenvalue of ${\cal L}_{\phi }$ is a
complex number $\lambda$
such that ${\cal L}_{\phi }v=\lambda v$ for a nonzero vector field $v \in
{\cal B}$. D. Ruelle [R1] proved the following theorem:

\vskip5pt
{\bf Theorem A.} {\em
The operator ${\cal L}_{\phi }$ has a positive and single maximal eigenvalue $\lambda
$ with an eigenvector $h$, and the remainder of the spectra
 of
${\cal L}_{\phi }$ is contained in a disk of radius strictly
less than $\lambda$.
Moreover, if $\sigma $, $\phi $ are $C^{k}$ for $k = 1$, $2$, $\cdots $,
$\omega$, then $h$ is a $C^{k}$ vector field.}

\vskip5pt
More recently, in [P],
[T] and [R2], the function $\phi $ was allowed
to be a complex function and
the spectral radius and
the essential spectral radius of ${\cal L}_{\phi }$ on the
$C^{k+\alpha }-$setting for $k=1$, $2$, $\cdots $, $0\leq \alpha \leq
1$ were estimated.

\vskip5pt
{\bf Feigenbaum's universality.}
Consider a family of unimodal mappings defined on
$[-1,1]$,
which is like the family $f_{t}(x) = t-(1+t)x^{2}$ for $t\in [0,1]$.
Suppose $t_{n}$ is the bifurcation value of parameters $t$ such that $f_{t}$,
$t<t_{n}$, does not have any periodic orbit of period $2^{n}$
and $f_{t}$, $t>t_{n}$, has an periodic orbit of period $2^{n}$.
M. Feigenbaum [F], and independently,  P. Coullet and C. Tresser [CT],
observed that the ratio
$\delta_{n}= \frac{t_{n}-t_{n-1}}{t_{n+1}-t_{n}}$ converges to a universal
number $\delta=4.669\cdots $ as $n$ goes to infinity. To explain this
universality, Feigenbaum [F] posed the following conjecture.

\vskip5pt
Suppose $f:[-1,1] \rightarrow [-1,1]$ is a symmetric
analytic folding mapping with a unique non-degenerate critical point $0$ and
satisfies $f(0)=1$ and
$f^{\circ 2}(0) <0 < f^{\circ 4}(0) <-f^{\circ 2}(0) < f^{\circ
3}(0) $. Let $q=f^{\circ 2}(0)$ and $I_{q}=[q, -q]$. The mapping
$f\circ f|I_{q}: I_{q}\rightarrow I_{q}$ is again a folding mapping
with a unique non-degenerate critical point.  Suppose $\alpha_{f}$ is the
 linear
rescaling of $I_{q}$ to $[-1,1]$ with $\alpha_{f}(q)=-1$.  Then
$F=\alpha_{f}\circ f\circ
f\circ \alpha_{f}^{-1}$ is a symmetric analytic folding mapping defined
on $[-1, 1]$. Denote $F$ by $R(f)$. Then $R$ is called the period doubling operator.

\vskip5pt
{\bf Conjecture A.} {\em The operator $R$ has a hyperbolic fixed point $g$
with (i) codimension one contracting manifold and (ii) dimension one
expanding manifold.}

\vskip5pt
O. Lanford proved this conjecture
with computer
assistance [L1]. After
him some mathematicians proved the existence of the fixed
point of $R$ without computer assistance,
for example, [CE] and [E].  Recently,
one of us proved the existence of the fixed point
$g$ and part (i) using quasi-conformal theory [S]. The proof in [S] not only
works for the period doubling operator but also works for its generalization, the renormalization operator (see Remarks in this introduction).
However, part $(ii)$ still lacks of a conceptual proof
(which hopefully, is valid also for the periodic
points of the renormalization operator). 
O. Lanford [L2] asked for a completely conceptual proof.

\vskip5pt
{\bf What we would like to say in this paper.}
We
give a proof of part $(ii)$. We use an induced operator
${\cal
L}_{\varphi }$ to study the linearization $T_{g}R$ of the period doubling
operator $R$ at the fixed point $g$. The operator ${\cal L}_{\varphi }$ is a ``Perron-Frobenius type operator'', but it is not a positive operator. The eigenvalues of
${\cal L}_{\varphi }$ agree with the eigenvalues of the linearization
$T_{g}R$ except for the value $1$. We use the
linear operator ${\cal L}_{n}$ with the finite rank $2^{n-1}$ to approximate
${\cal L}_{\varphi }$ in the $C^{b}$-setting ($C^{b}$ is the space of
bounded vector fields on $g(I)$).
Under {\em the assumption that
$g$ is a concave function [L1]},
we prove the following statements:

\vskip3pt
$(1)$ Each ${\cal L}_{n}$ has
an
eigenvalue $\lambda_{n} >1$ with a positive eigenvector $v_{n}$, this means
that each component of $v_{n}$ is positive.

\vskip3pt
$(2)$
There is a subsequence $\{ n_{i} \}_{i=0}^{\infty}$ of the integers
such that the limit $\lambda =\lim_{i\rightarrow
\infty} \lambda_{n_{i}}>1$ is an eigenvalue of ${\cal L}_{\varphi }$
with an eigenvector $v=\lim_{n\rightarrow +\infty }v_{n_{i}}$ in
$C^{b}$.

\vskip3pt
$(3)$ The number $\lambda $ is an eigenvalue of ${\cal L}_{\varphi }$ in
the $C^{0,1}$-setting ($C^{0,1}$ is the space of Lipschitz continuous
vector fields on $g(I)$ ).

\vskip3pt
$(4)$ The limit
$\lambda $ is an eigenvalue of ${\cal L}_{\varphi }$ in
the $C^{\omega}$-setting ($C^{\omega}$ is the space of analytic vector fields on
$g(I)$).

\vskip5pt
These yield a proof of part $(ii)$. The proofs are constructive.  One can
calculate the
approximating expanding manifolds and the rate of expansion of $R$ by using
${\cal L}_{n}$.

\vskip5pt
We
also learned that recently, J.-P. Eckmann and H. Epstein [EE] gave a
different proof of part $(ii)$
and R. Artuso, E. Aurell and P. Cvitanovi\`c [AAC] gave a
rigorous mathematical proof of part $(ii)$.

\vskip5pt
{\bf Some remarks on the renormalization operator.}
Suppose $f:[-1,1] \rightarrow [-1,1]$ is a symmetric
analytic unimodal mapping with a unique non-degenerate critical point $0$.
Suppose there is an integer $n>1$ such that there exists an interval
$I$ containing $0$ and the restriction of the $n$-fold $f^{\circ n}$ of $f$
 maps $I$ into itself.	Let $n$ be the smallest such integer
and $I_{q}=[q, -q]$ or $[-q, q]$ be the
maximal such interval. The point $q$ is a fixed point of $f^{\circ n}$.
Let $\alpha_{f}$ be the linear
mapping which rescales $I_{q}$ to $[-1,1]$ with $\alpha_{f}(q)=-1$.  Then
$F=\alpha_{f}\circ f^{\circ n}\circ \alpha_{f}^{-1}$ is a symmetric
analytic unimodal mapping defined
on $[-1, 1]$.  We say $f$ is once renormalizable and
$R: f\mapsto F$ is the renormalization operator.

\vskip5pt
{\bf Conjecture B.} {\em $(I)$ For every periodic kneading sequence $\rho
=(w_{1}*w_{2}*\cdots *w_{k})^{*\infty}$, where $\rho$ is decomposed into the
star product of primary sequences, (see [MT] and [CEc] for a definition of 
a kneading sequence, a definition of star product and a definition of a
primary sequence), there is a hyperbolic periodic point $g_{\rho }$
(with this kneading sequence)
of period $k$ of $R$ with (i) codimension one contracting manifold and $(ii)$
dimension one expanding manifold.

Moreover, $(II)$ $R$ is hyperbolic on its maximal invariant set with
$(i)$ codimension one stable manifold and $(ii)$ dimension one unstable
manifold.
}

\vskip5pt
Topologically, we knew that the maximal invariant set of
$R$ is like the Smale horse shoe. Under
the assumption that $g_{\rho}$ is a concave function and some a prior
estimate of linear rescale mapping $\alpha_{g_{\rho}}$,
one may use the methods in this paper to prove part $(ii)$ of $(I)$.  But as
H. Epstein pointed out to us if we also consider a power law critical
point and the exponent of $g$ at its power law critical point is large, then $g$
is not a concave function any more.
In this paper, the concave condition is used only in the proof of
statement $(1)$.
We are expecting a proof of statement $(1)$ without
the assumption that $g_{\rho}$ is a concave function.	This seems to be a
promising problem.
The other option is to prove that $g_{\rho}$
is a concave function for every periodic kneading sequence
$\rho$ in the case that the critical point of $g_{\rho}$ is
non-degenerate. But it seems to be a difficult problem.

\vskip3pt
The other observation is that the generalized Feigenbaum's
$\delta_{\rho }$ only depends on the
data related to the critical orbit.

\vskip20pt
\noindent {\bf Acknowledgement:}
  We
would like to thank Folkert Tangerman for many helpful conversations and
Viviane Baladi for many useful suggestions and comments to this paper.
We would also like to thank John Milnor for reading this paper and for his
many helpful remarks. 
This manuscript was first written on June, 1989 when one of
us visited IHES.
He would like to thank IHES for its hospitality.

\vskip10pt
\centerline{\large \bf \S 2 The Period Doubling Operator and the Induced Operator}

\vskip5pt
Suppose $I$ is the interval $[-1, 1]$ and $U\subset {\bf C}^{1}$ is a
connected open subset
containing $I$.  Let ${\cal B}(I, U)$ be the space of
folding mappings $f$ from $I$ into $I$ with a unique non-degenerate critical
point $0$ and an analytic extension $F$ on $U$ which can be extended
to the boundary $\partial U$ continuously.  Suppose
${\cal B}_{s}(I, U)$ is the subspace of even functions in ${\cal
B}(I, U)$ and ${\cal B}_{s,0}(I, U)$ is the subspace of mappings which
are in ${\cal B}_{s}(I, U)$ and satisfy the conditions
$ f(0)=1$ and $f^{\circ
3}(0) >-f^{\circ 2}(0)> f^{\circ 4}(0) >0 > f^{\circ 2}(0)$. The period
doubling operator $R$ from ${\cal B}_{s,0}(I, U)$ into ${\cal
B}_{s}(I, U)$ is defined by
\[ R(f)(x) = -\alpha_{f}f\circ f(-\alpha_{f}^{-1}x), \hskip8pt x\in
I,\] for $f\in {\cal B}_{s,0}(I,U)$, where $\alpha_{f}=-\frac{1}{f(1)}$.

\vskip5pt
Suppose $g$ is the fixed point of $R$ ([VSK]) and $U$
is an open
set contained in the domain of $g$.  Let $\alpha =-\frac{1}{g(1)}$, $J=
g(I)$ and $\Omega = g(U)$.  Suppose
${\cal V}^{\omega}(J, \Omega )$ is the space of real vector fields $v$ on $J$ with
a complex analytic extension $V$ on $\Omega$ which can be extended to the
boundary $\partial \Omega$ continuously.
This space equipped with the uniformly
convergent norm is a Banach space.

\vskip5pt
\noindent {\bf \S 2.1 From the period doubling operator to the induced
operator}

\vskip5pt
Suppose $J_{0}$ and $J_{1}$ are the intervals $[g(1), g^{\circ 3}(1)]$
and $[g^{\circ 2}(1),1]$. We define $\sigma$ from 
$J_{0}\cup J_{1}$ onto $J$ by
\[ \sigma (x) =\left\{ \begin{array}{ll}
			 -\alpha x ,& x\in J_{0},\\
			 -\alpha g(x) ,& x\in J_{1}.
			 \end{array}
		\right.\]
The mapping $\sigma$ is
expanding with expansion constant $\alpha $ for
$|g^{'}(x)| >1$, $x\in J_{1}$, and it has an analytic extension, which we
still denote as $\sigma$, on
$\Omega_{0}\cup \Omega_{1} \supset J_{0}\cup
J_{1}$ with also expansion constant $\alpha $.
Here $\Omega_{0}$ and $\Omega_{1}$ are disjoint subdomains of $\Omega$ and contain
$J_{0}$
and $J_{1}$, respectively. Moreover, the restrictions $\sigma | \Omega_{0}$ and
$\sigma |\Omega_{1}$ of $\sigma$ to $\Omega_{0}$ and $\Omega_{1}$ are
bijective from $\Omega_{0}$ and $\Omega_{1}$ to $\Omega$ and can be extended 
continuously to the boundaries $\partial \Omega_{0}$ and $\partial \Omega_{1}$, respectively (see Figure 1).

\centerline{\psfig{figure=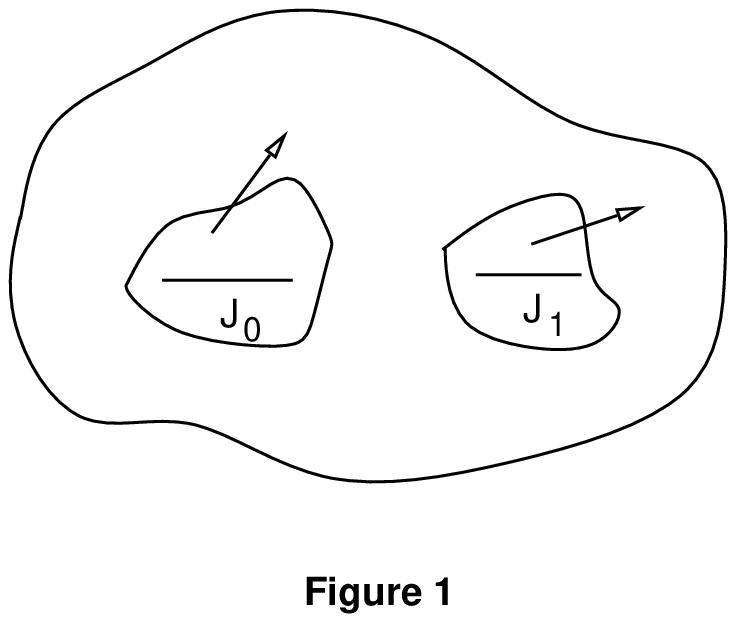}}

Suppose $C$ is the attractor of $g$ and $\Lambda$ is the maximal
invariant set of $\sigma$.

\vskip5pt
{\bf Lemma 1.} {\em The set $\Lambda$ and the set $C$ are the same.}

\vskip5pt
{\it Proof.} The reader may check it by the equation $g(x)=-\alpha
g\circ g(-\alpha^{-1}x)$.

\vskip5pt
Suppose $\varphi(z)$ is the derivative $\sigma^{'}(z)$ of $\sigma $
on $\Omega_{0}\cup \Omega_{1}$. We define
${\cal L}_{\varphi }$ from ${\cal V}^{\omega}(J, \Omega )$ into ${\cal
V}^{\omega}(J, \Omega )$ by
\[ ({\cal L}_{\varphi }v)(z) =\sum_{w\in \sigma^{-1}(z)} \varphi
(w)v(w)\]
and call it the induced operator.  It is a ``Perron-Frobenius type
operator'' but is not positive.  It is clearly bounded and
compact (by Montel's theorem).  

\vskip5pt
Suppose $T_{g}{\cal B}_{s,0}(I, U)$
is the tangent space of ${\cal B}_{s,0}(I, U)$ at $g$ and
$T_{g}R$ from $T_{g}{\cal B}_{s,0}(I,
U)$ into $T_{g}{\cal B}_{s,0}(I, U)$ $(=T_{g}{\cal B}_{s}(I, U))$ is the
tangent map of $R$ at $g$.

\vskip5pt
{\bf Lemma 2.} {\em The mapping $g_{*}$ from ${\cal V}^{\omega}(J, \Omega )$
into $T_{g}{\cal B}_{s,0}(I, U)$ defined by
$g_{*}(v)(x) = v(g(x))$ for $x\in \Omega $ and $v\in {\cal
V}^{\omega}(J, \Omega)$
is an isomorphism.}

\vskip5pt
{\it Proof.} The proof is easy.

\vskip5pt
{\bf Lemma 3.} {\em The operators
${\cal L}_{\varphi }$ and $T_{g}R$ have the same
eigenvalues (counted with multiplicity) except for the value $1$.}

\vskip5pt
{\it Proof.} By some calculations, we can show that
\[ {\cal L}_{\varphi }=g_{*}^{-1}\circ T_{g}R\circ g_{*}+e_{1},\]
where $e_{1}$ is the projection from ${\cal V}^{\omega}(J, \Omega)$ to the
eigenspace of eigenvalue one.

\vskip5pt
{\it Remark.}
Suppose $V_{2m-1}(x)= g^{'}(x)x^{2m-1} -(g(x))^{2m-1}\in T_{g}{\cal
B}_{s,0}(I, U)$ and $v_{2m-1}=g_{*}^{-1}(V_{2m-1})\in {\cal V}^{\omega}(J,
\Omega)$. The vector
$v_{2m-1}$ is an eigenvector of ${\cal L}_{\varphi }$ with eigenvalue
$\lambda_{2m-1} =\alpha^{-(2m-2)}$ for $ m=1, 2,\cdots $.

\vskip5pt
Lemma 3 tells us that we can use ${\cal L}_{\varphi }$ which has an 
explicit form to study the eigenvectors and eigenvalues of $T_{g}R$ except the value $1$.
We will use it to find the expanding direction and the rate of $R$
at the fixed point $g$.

\vskip15pt
\noindent {\bf \S 2.2 The induced operator ${\cal L}_{\varphi}$}

\vskip5pt
Suppose $v$ is a real vector field on $\Lambda $. We say it is a Lipschitz
continuous if there is a constant $M>0$ such that $|v(x)-v(y)| \leq
M|x-y|$ for any $x$ and $y$ in $\Lambda $. We say it is bounded if there
is a constant $M>0$ such that $|v(x)|\leq M$ for any $x$ in $\Lambda$.
Let ${\cal V}^{0,1}(\Lambda )$ be the space of real Lipschitz continuous
vector fields on $\Lambda $ and ${\cal V}^{b}(\Lambda )$ be the space of
bounded vector fields on $\Lambda $.  Suppose $\varphi (x)$ is the
derivative $\sigma^{'}(x)$ on $\Lambda $. We define two linear operator
by the same formula. One is ${\cal L}_{\varphi , L}$ from 
${\cal V}^{0,1}(\Lambda )$ into ${\cal V}^{0,1}(\Lambda )$ defined by
\[
({\cal L}_{\varphi, L}v)(x) = \sum_{y\in \sigma^{-1}(x)} \varphi
(y)v(y)\]
and the other is
${\cal L}_{\varphi , B}$ from 
${\cal V}^{b}(\Lambda )$ into ${\cal V}^{b}(\Lambda )$ defined by
\[
({\cal L}_{\varphi, B}v)(x) = \sum_{y\in \sigma^{-1}(x)} \varphi
(y)v(y).\]
They are bounded but not compact.

\vskip5pt
{\bf Lemma 4.} {\em Suppose $\lambda $ is an eigenvalue of
${\cal L}_{\varphi , B}$ and $\lambda > \alpha +1$. Then it is an eigenvalue of
${\cal L}_{\varphi , L}$.}

\vskip5pt
{\it Proof.} There is a nonzero vector field $v$ in ${\cal
V}^{b}(\Lambda )$ such that
\[ {\cal L}_{\varphi , B}v=\lambda v.\]
This is
\[ -\alpha v(-\alpha^{-1}x)-\alpha g'(g^{-1}(-\alpha^{-1}x))
v(g^{-1}(-\alpha^{-1}x))=\lambda v(x) \hskip 25pt (*) \]
for any $x$ in $\Lambda $.  From this we can have an inequality
\[ \max_{x\neq y\in \Lambda }(\frac{|v(x)-v(y)|}{|x-y|}) \leq
\frac{M}{\lambda -\alpha -1}\]
where $M$ is a positive constant. In the other words, $v$ is Lipschitz
continuous on $\Lambda $ and is an eigenvector of
${\cal L}_{\varphi , L}$
with the eigenvalue $\lambda $.

\vskip5pt
{\bf Lemma 5.} {\em Suppose $\lambda $ is an eigenvalue of
${\cal L}_{\varphi , L}$ and $\lambda >\alpha +1$. Then $\lambda $ is an
eigenvalue of ${\cal L}_{\varphi }$.}

\vskip5pt
{\it Proof.} The basic idea to proof this lemma is to use the fact that
the grand preimage $\cup_{n=0}^{\infty}g^{-n}(\Lambda )$ is a dense
subset on $I$ and to use the equality (*) countably many times.
We will not write down our proof in detail because recently, there is a
more general
theorem proved by D. Ruelle [R2]. One of us learned this theorem
when he visited IHES. We outline some Ruelle's result here.

\vskip5pt
\noindent {\bf \S 2.3 A general theorem for operators like the
induced operator}

\vskip5pt
In this subsection, the notations $J_{0}$, $J_{1}$, $\Omega$, $\Omega_{0}$
and $\Omega_{1}$ are the same as that in \S 2.1. 

\vskip5pt
Suppose $e$ from $J_{0}\cup J_{1}$ into and onto $J$ is an expanding mapping
such that the restrictions $e|J_{0}$ and $e|J_{1}$ of $e$ to $J_{0}$ and
$J_{1}$ have bijective,
expanding, analytic extensions $F_{0}$ and $F_{1}$ from $\Omega_{0}$ and $\Omega_{1}$ to
$\Omega$, respectively. Moreover, $F_{0}$ and $F_{1}$ can be extended 
to the boundaries
$\partial \Omega_{0}$ and  
$\partial \Omega_{1}$ continuously, We use $E$ to denote the expanding map
\[ E(z) =\left\{ \begin{array}{ll}
			 F_{0}(z) ,& z\in \Omega_{0},\\
			 F_{1}(z) ,& z\in \Omega_{1}
			 \end{array}
		\right. \]
from $\Omega_{0}\cup \Omega_{1}$ to $\Omega$.

\vskip5pt
Suppose $\phi$ from $J_{0}\cup
J_{1}$ into the real line is a real analytic function with a complex
analytic extension $\Phi$ from $\Omega_{0}\cup \Omega_{1}$ into the complex
plane which
can be also extended continuously to the boundaries $\partial
\Omega_{0}\cup \partial \Omega_{1}$.
Let $\theta$ be the expanding
constant of $E$ and $\Lambda_{E}$ be its maximal invariant set.
A linear operator ${\cal
L}_{\phi }$ from ${\cal V}^{\omega}(J, \Omega )$ into ${\cal V}^{\omega}(J, \Omega )$
is defined by
\[ ({\cal L}_{\phi})(v)(z) =\sum_{w\in E^{-1}(z)}\phi (z)v(z).\]
Suppose ${\cal V}^{0,1}(\Lambda_{E} )$ is the space of real Lipschitz continuous
vector fields on $\Lambda_{E}$.
Let $|\phi|$ be the function which takes values $|\phi (z)|$ at all
$z\in \Omega_{0}\cup \Omega_{1} $.  We assume $|\phi |$ is a positive function. Then
the operator ${\cal L}_{|\phi |}$ from ${\cal V}^{0,1}(\Lambda_{E} )$ into
${\cal
V}^{0,1}(\Lambda_{E} )$ defined by
\[
({\cal L}_{|\phi|})(v)(x) =\sum_{y\in E^{-1}(x)}|\phi |(y)v(y)\]
is an
Perron-Frobenius operator.
For the positive function $|\phi |$, we can define its pressure as
\[ P(\log |\phi|) = \sup_{\mu}(h_{\mu}(E) +\int_{\Lambda_{E}}(\log |\phi
|)d\mu )\]
where $\mu $ is an invariant measure of $E$ and $h_{\mu }$ is the
measure-theoretic entropy of $E$ with respect to $\mu $.
By the variation principle (see, for example, [B]), we have that
\[ P(\log |\phi
|)
=\lim_{n\rightarrow \infty }\frac{1}{n}\log (\sum_{x\in
fix(E^{\circ n})}\prod_{i=0}^{n-1} |\phi (E^{\circ i}(x))|).\]
Let $A=exp(P(\log |\phi |))$. It is a simple eigenvalue of
${\cal L}_{|\phi|}$ and
all other eigenvalues of ${\cal L}_{|\phi|}$ are in the open disk
$D_{A}=\{ z|\in {\bf C}^{1}, |z| < A \}$ (see the theorem in the
introduction). Suppose ${\cal
L}_{\phi ,L}$ from ${\cal V}^{0,1}(\Lambda_{E})$ into ${\cal V}^{0,1}(\Lambda_{E})
$ is defined by
\[ ({\cal L}_{\phi , L})(v)(x) =\sum_{y\in E^{-1}(x)}\phi (y)v(y)\]
and $A_{1}=\theta^{-1}A$.

\vskip5pt
{\bf Lemma 6.} (see [R2]). {\em All
the eigenvalues of ${\cal L}_{\phi , L}$ are in the open disk $D_{A}$
and if $\lambda $ is an eigenvalue
of ${\cal L}_{\phi , L}$ and $|\lambda |> A_{1}$, then $\lambda $ is an
eigenvalue of ${\cal L}_{\phi }$.}

\vskip5pt
Suppose $\sigma$ is the expanding
mapping induced from the period doubling operator and
$\varphi$ is the derivative $\sigma^{'}$. The expanding constant of
$\sigma$
is $\alpha $. By some combinatorial arguments, we have that 
\[ \sum_{x\in fix(\sigma^{\circ n})} \prod_{i=0}^{n-1} \log |\phi
(\sigma^{\circ i}(x))| \leq ( \alpha^{2} +\alpha )^{n}\]
Moreover, 
by using the variation principle,
\[ A=
exp(P(\log |\varphi
|)) \leq \alpha (\alpha +1)\]
and thus $A_{1} = \alpha^{-1}A
\leq \alpha +1$. From this, Lemma 6 gives a proof of Lemma 5.
Moreover, if $\lambda $ is an
eigenvalue of ${\cal L}_{\phi , B}$ and $\lambda > \alpha +1$, then it
is an eigenvalue of $T_{g}R$.

\vskip10pt
\centerline{\Large \bf \S 3 The Construction of the Expanding Direction}

\vskip5pt
We prove that the induced operator ${\cal L}_{\varphi}$ has an
expanding direction and construct this direction in this section. The transformation of
this direction under $g_{*}$ is the expanding direction of the period
doubling operator.

\vskip5pt
\noindent {\bf \S 3.1 An easy observation}

\vskip5pt
Suppose $I=[-1, 1]$ is a closed interval of the real line ${\bf R}^{1}$ and
$D \supset I$ is
an open disk in the complex plan
${\bf C}^{1}$. Suppose $I_{0}$ and $I_{1}$ are disjoint closed 
subintervals of $I$ and $e$ is a piecewise linear expanding map from
$I_{0}\cup I_{1}$ onto and into $I$ with the derivative
\[ \phi (x) =e'(x) =\left\{
\begin{array}{ll} 				-a, & x\in I_{0},\\
				 b, & x\in I_{1}
				\end{array}
		     \right.\]
where $b>a >2$ are two constants.
Let $E$ be the extension of $e$ from $D_{0}$ and $D_{1}$ onto and into $D$
with also the derivative, we still denote it by $\phi$,
\[ \phi (z) =E'(z) =\left\{
\begin{array}{ll} 				-a, & z\in D_{0},\\
				 b, & z\in D_{1}.
				\end{array}
		     \right.\]

\vskip5pt
Let
${\cal V}^{\omega}(I, D)$ be the space of real vector fields $v$ on $I$ with
a complex analytic extension $V$ on $D$ which can be extended to the
boundary $\partial D$ continuously.
For $t\in [0,1]$, we define
\[ \phi_{t} (x) =\left\{ \begin{array}{ll}
				ae^{2\pi it}, & x\in J_{0},\\
				 b, & x\in J_{1}
				\end{array}
		     \right.\]
and ${\cal L}_{\phi_{t}}$ from ${\cal V}^{\omega}(I, D )$ into ${\cal
V}^{\omega}(I, D )$ by
\[ ({\cal L}_{\phi_{t}})(v)(z) =ae^{2\pi it } v(E_{0}(z))
+bv(E_{1}(z))\]
where $E_{0}$ and $E_{1}$ are the
inverse branches of $E$.

\vskip5pt
{\bf Proposition A.}
{\em The set $\{ \lambda_{n,t}= \frac{1}{b^{n-1}}+
\frac{1}{e^{2\pi it(n-1)}a^{n-1}}\}_{n=0}^{\infty }$ is the spectrum of
${\cal L}_{\phi_{t}}$.}

\vskip5pt
{\sl Proof.} Suppose the center of $D$ is $0$. Then every $v\in {\cal
V}^{\omega}(I,D)$ has the Taylor expansion 
\[ v(z)=\sum_{k=0}^{\infty} a_{k} z^{k},\]
where $a_{k}$ are all real numbers. To find $\lambda_{n,t}$ for $n=0$,
$\cdots $, $+\infty$, we may solve
the equation
\[ {\cal L}_{\phi_{t}} v_{n} =\lambda_{n,t}v_{n}\]
for $v_{n}(z)=\sum_{k=0}^{n}a_{k}z^{k} \in {\cal V}^{\omega}(I, D)$ . 

\vskip5pt
Under the condition $b>a >2$,  $||\lambda_{n,t}|| <1$ for $n >1$.
The other two eigenvalues $\lambda_{1,t} (=2$ for all
$t\in [0,1]$) and
$\lambda_{0,t}= b+e^{2\pi it}a$ are special.  Here $2=\exp (h(E))$ is a
topological invariant where $h(E)$ is the topological
entropy of $E$.

\vskip5pt
From this proposition, we can observe that $\lambda_{0,t}$ is the 
maximal eigenvalue of
${\cal L}_{\phi_{t}}$ for all $t\in [0,1]$
if and only if $b-a >2$.
In the other words, It is the maximal eigenvalue of
${\cal L}_{\phi_{t}}$ 
for all
$t\in [0,1]$ if and only if $b-a >2$: an unbalanced condition, the
orientation preserving part is much stronger than the orientation
reversing part or vice versa.

\vskip5pt
\noindent {\bf \S 3.2 The construction}

\vskip5pt
We use the same notations as that in the previous section. We
note that the derivative of $\sigma$ on $J_{1}$ is
strictly greater than one and the derivative of $\sigma$
at the right end point of $J$ is $\alpha^{2}$.

\vskip5pt
Suppose $\varphi_{1}$ is the function defined by
\[ \varphi_{1}=\left\{ \begin{array}{ll}
			 -\alpha  , & x\in J_{0}\cap \Lambda,\\
			\alpha^{2}  , & x\in J_{1}\cap \Lambda
			\end{array}
		\right.\]
and ${\cal L}_{1}:{\cal V}^{b}(\Lambda )\rightarrow {\cal V}^{b}(\Lambda
)$ is the corresponding operator defined by
\[ ({\cal L}_{1}v)(x)= \sum_{y\in \sigma^{-1}(x)}\varphi_{1}(y)v(y).\]
The number $\lambda_{1}
=\alpha (\alpha -1)$ is an eigenvalue of ${\cal
L}_{1}$ with an eigenvector $v_{1}=1$ on $\Lambda $.

\vskip5pt
Suppose
$\sigma^{-2}(J) =J_{21}\cup J_{22}\cup J_{23}\cup J_{24}$ and
$J_{23}=[a_{21}, b_{21}]$, $J_{24}=[a_{22}, b_{22}]$ (see Figure 1).
Let
$\beta_{21}=|g^{\prime }(b_{21})|$ and
$\beta_{22}=|g^{\prime }(b_{22})|=|g^{\prime }(1)|=\alpha $.  Because $g$ is a concave function [L1], we have that $\beta_{21}\leq
\beta_{22}$.  Suppose $\varphi_{2}$ is the function defined by
\[ \varphi_{2}(x)=\left\{ \begin{array}{ll}
			 -\alpha , & x\in J_{0}\cap \Lambda,\\
			\alpha \beta_{21} , & x\in J_{23}\cap \Lambda ,\\
			\alpha \beta_{22} , & x\in J_{24}\cap \Lambda
			\end{array}
		\right.\]
and ${\cal L}_{2}:{\cal V}^{b}(\Lambda )\rightarrow {\cal V}^{b}(\Lambda
)$ is the corresponding operator defined by
\[ ({\cal L}_{2}v)(x)= \sum_{y\in \sigma^{-1}(x)}\varphi_{2}(y)v(y).\]
Let $k_{21}$ be the vector field on $\Lambda $ defined by
\[ k_{21}(x)= \left\{ \begin{array}{ll}
			    1 ,& x\in J_{0}\cap \Lambda,\\
			    0 ,& x\in J_{1}\cap \Lambda
			   \end{array}
		   \right.\]
and $k_{22} =1-k_{21}$. The space ${\bf R}^{2} =span \{ k_{21},
k_{22}\} $ is a subspace of ${\cal V}^{b}(\Lambda )$.
For any $v = x_{21}k_{21} +x_{22}k_{22}$,
\[ ({\cal L}_{2}v)(x)= (k_{21}, k_{22}) \left(
\begin{array}{ll}
      -\alpha , & \alpha \beta_{21} \\
      -\alpha , & \alpha \beta_{22}
\end{array}
\right)
\left(
\begin{array}{l}
      x_{21} \\
      x_{22}
\end{array}
\right).\]
Let $A_{2}$ be the matrix
\[ \left(
\begin{array}{ll}
      -\alpha , & \alpha \beta_{21} \\
      -\alpha , & \alpha \beta_{22}
\end{array}
\right). \]

\vskip5pt
{\bf Proposition C.} {\em The maximal 
eigenvalue of $A_{2}$ is
\[ \lambda_{2} =\alpha \frac{(\beta_{22}-1) +\sqrt{(\beta_{22}-1)^{2}
+4(\beta_{22}-\beta_{21})}}{2},\]
with an eigenvector $v_{2}=(t_{21}, 1)$, $t_{21}<1$.}

\vskip5pt
{\it Proof.} The proof uses linear algebra.

\vskip5pt
Furthermore,
suppose $\sigma^{-n}(J) =J_{n1}\cup
J_{n2}\cup \cdots \cup J_{n2^{n-1}}\cup J_{n(2^{n-1}+1)}\cup \cdots \cup J_{n 2^{n}}$ and
$J_{n(2^{n-1}+i)} =[a_{ni}, b_{ni}]$ (see Figure 2).

\vskip5pt
\centerline{\psfig{figure=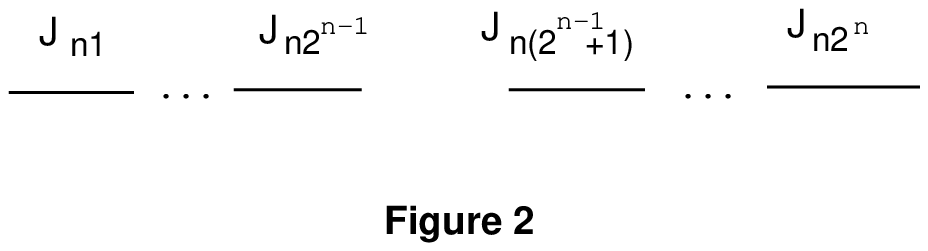}}

\vskip5pt
Let $\beta_{ni}
=|g^{\prime}(b_{ni})|$ for $i=1, 2, \cdots 2^{n-1}$.  Because $g$ is a concave function, we have that
\[ 1< \beta_{n1} <\cdots <\beta_{n 2^{n-1}} =\alpha .\]
Suppose $\varphi_{n}$ is the function defined by
\[ \varphi_{n}(x)=\left\{
\begin{array}{ll}
      -\alpha , & x \in J_{0}\cap \Lambda,\\
      -\alpha \beta_{ni}  ,& x\in J_{n(2^{n-1}+i)}\cap \Lambda ,\hskip7pt
i=1,2,\cdots 2^{n-1},
 \end{array}
\right.\]
and ${\cal L}_{n}$ from ${\cal V}^{b}(\Lambda )$ into ${\cal V}^{b}(\Lambda
)$ is the corresponding operator defined by
\[ ({\cal L}_{n}v)(x)= \sum_{y\in \sigma^{-1}(x)}\varphi_{n}(y)v(y).\]
Let $k_{ni}$ be the vector field on $\Lambda $ defined by
\[ k_{ni}(x)= \left\{ \begin{array}{ll}
			    1 , & x\in (J_{n(2i-1)}\cup J_{n(2i)})\cap
\Lambda,\\
			    0 , & x \in \Lambda\setminus
((J_{n(2i-1)}\cup J_{n(2i)})\cap
\Lambda)
			   \end{array}
		   \right.\]
for $i= 1$, $2,$ $\cdots $, $2^{n-1}$. The space ${\bf R}^{2^{n-1}} =span \{
k_{n1},\cdots , k_{n2^{n-1}}\}$ is a subspace of ${\cal V}^{b}(\Lambda )$.
For any
$v=x_{n1}k_{n1}+\cdots x_{n2^{n-1}}k_{n2^{n-1}} $, we have that
\[ {\cal L}_{n}v = K_{n} A_{n} X_{n}^{t}\]
where $K_{n}=(k_{n1}, \cdots , k_{n2^{n-1}})$ and $X_{n}=(x_{n1},
\cdots , x_{n2^{n}})$ and $A_{n} $ stands for the
$2^{n-1}\times 2^{n-1}$-matrix \[ \left( \begin{array}{llllllllll}
		0 & 0 & \cdots & 0 & -\alpha & \alpha \beta_{n1} &
0 &\cdots &0 & 0\\
		0 & 0 & \cdots  & 0  & -\alpha  & \alpha \beta_{n2}  &
0 & \cdots  & 0  & 0\\
		0 & 0 & \cdots  & -\alpha  & 0  & 0  & \alpha \beta_{n3}
& \cdots  & 0  & 0\\
		0 & 0 & \cdots  & -\alpha  & 0  & 0  & \alpha \beta_{n4} &
\cdots	& 0  & 0\\
		\vdots & \vdots & \cdots & \vdots  & \vdots  & \vdots
& \vdots  & \cdots  & \vdots & \vdots\\
		0 & -\alpha & \cdots	& 0  & 0  & 0
 & 0  &\cdots  &\alpha \beta_{n(2^{n-1}-3)}  & 0\\
		0 & -\alpha & \cdots & 0  & 0  & 0
  & 0 & \cdots &\alpha \beta_{n(2^{n-1}-2)}  & 0\\
		-\alpha  & 0 & \cdots  & 0  &0  & 0  &
  & \cdots & 0	& \alpha \beta_{n(2^{n-1}-1)}\\
		-\alpha & 0 & \cdots	& 0  & 0  & 0
  & 0  &\cdots & 0  &\alpha \beta_{n2^{n-1}}
\end{array}
\right).\]

\vskip5pt
{\bf Proposition D.}
{\em The matrix $A_{n}$ has an eigenvalue $\lambda_{n}$
which is greater than $\alpha (\alpha -1)$.}

\vskip5pt
{\it Proof.} Suppose $CN_{n}$ is the set $\{ (x_{n1}, \cdots ,
x_{n2^{n-1}})| \in {\bf R}^{2^{n-1}},$
$x_{ni}\geq 0$ for $i=1$, $\cdots $, $i=2^{n-1}$
and $x_{n1}\leq x_{n2}\leq \cdots \leq x_{n2^{n-1}}\}$. It is easy to
check that $CN_{n}$ is a convex cone and $A_{n}$ maps this cone into
the interior of this cone and zero vector. By the Brouwer fixed point
theorem, we conclude that there is a unique direction
${\bf R}^{+}v_{n}$ in this cone which is preserved by
$A_{n}$.
Suppose $v_{n}=(t_{n1}, \cdots , t_{n2^{n-1}})$
with $t_{n2^{n-1}}=1$ is an eigenvector with the eigenvalue
$\lambda_{n}$.
By the equation $A_{n}v_{n} =\lambda_{n}v_{n}$, we have
that $-\alpha t_{n1}+\alpha^{2}=\lambda_{
n}$. Because $t_{n1}<1$, we get $\lambda_{n} >\alpha (\alpha -1).$

\vskip5pt
{\it Remark.} Because the cone $CN_{n}$ is a subset
of the cone $CN_{n+1}$ for any $n\geq 1$, we can prove more that $\{
\lambda_{n}\}_{n=1}^{\infty }$ is an increasing sequence. But we will
not use this fact because we would like the following arguments to be
true even in the case that $g$ is not a concave function.

\vskip5pt
{\bf Proposition E.} {\em There is a subsequence $\{ n_{i}\}$ of the
integers such that the continuous extension of the limit
$v=\lim_{i\mapsto \infty }v_{n_{i}}$ on the critical orbit $Or(g)=\{
g^{\circ
n}(0) \}_{n=1}^{\infty }$ is an eigenvector of ${\cal L}_{\varphi ,
B}$ with the eigenvalue $\lambda =\lim_{i\mapsto \infty
}\lambda_{n_{i}}$.}

\vskip5pt
{\it Proof.} Because $Or(g)$ is a countable set, we can find a
subsequence
$\{ n_{i}\}_{i=0}^{\infty}$ such that for every $a\in Or(g)$, the limit
$v_{n_{i}}(a)$ exists as $i$ goes to infinity. We denote this limit as
$v(a)$.  For the sequence $\{ \lambda_{n_{i}}\}_{i=0}^{\infty }$, we
can find convergent subsequence.
Let $\lambda$ be the limit of this subsequence. Then we have that
$({\cal L}_{\varphi ,B}v)(a)=\lambda v(a)$ for any $a\in Or(g)$. Now by
using the equation $(*)$ and the fact $\alpha (\alpha -1)> \alpha
+1$ which can be implied by $\alpha >1+\sqrt {2}$, we can show that $v$
has a continuous extension on $\Lambda $ which is the closure of
$Or(g)$.

\vskip5pt
\noindent {\bf \S 3.3 A program}

\vskip5pt
In \S 3.2, we use a subsequence of $\{ v_{n}\}_{n=0}^{\infty}$ to prove
that there is an expanding direction of ${\cal L}_{\varphi}$. Under the
assumption that $g$ is a concave function, we can say more on the
sequence
$\{ v_{n}\}_{n=0}^{\infty}$ and the corresponding eigenvalues $\{
\lambda_{n}\}_{n=0}^{\infty}$. For example, $\{ \lambda_{n} \}$ is an
increasing sequence and for every $a\in \Lambda $, $\{ v_{n}(a)\}$ is a
monotone sequence. In practice, we can use these good properties to
give an effective program to find the expanding
direction and the rate of the period doubling operator as follows:

\vskip5pt
Suppose $v$ is a vector in ${\bf R}^{k}$. We use $(v)_{i}$ to denote its
$i^{th}$-coordinate.

\vskip5pt
$(1)$ Start from the constant function $v_{1}=1$. Consider it as a
vector in ${\bf R}^{2}$ and compute the limiting
vector
\[ v_{2}=\lim_{l\mapsto \infty}
\frac{A_{2}^{l}v_{1}}{(A_{2}^{l}v_{1})_{2}}\]
and the corresponding eigenvalue $\lambda_{2}=\alpha (\alpha -
(A_{2}v_{2})_{1})$.

\vskip5pt
$(n)$ Let $v_{n-1}\in {\bf R}^{2^{n-2}}$ be the eigenvector of
$A_{n-1}$ with the eigenvalue $\lambda_{n-1}$ .
Consider $v_{n-1}$ as a vector in ${\bf R}^{2^{n-1}}$ and compute the
limiting vector
\[ v_{n}=\lim_{l\mapsto \infty}
\frac{A_{n}^{l}v_{n-1}}{(A_{n}^{l}v_{n-1})_{2^{n-1}}}\]
and the corresponding eigenvalue $\lambda_{n}=\alpha (\alpha -
(A_{n}v_{n})_{1})$.

$(\infty)$ The limiting vector
\[ V=\lim_{n\mapsto \infty}g_{*}(v_{n})\]
is the expanding direction and the limiting value
\[ \delta =\lim_{n\mapsto \infty}\lambda_{n} \]
is the rate of expansion of the period doubling operator at the fixed point
$g$.

\end{document}